\magnification\magstephalf
\input AHTOH-E.STY
\vsize24truecm
\hfuzz5.1pt
\def\Chi{{\rm X}}
\let\tmp\bigtimes
\def\bigtimes{\mathop{\!\!\tmp\!\!}}
\UDC{
512.543.72 
+
512.542.68 
+
512.544.33 
}

\title{
Finite and nilpotent strongly verbally closed groups
}

\author{%
Anton A. Klyachko$^{\flat\natural}$
\quad
Veronika Yu. Miroshnichenko$^\flat$
\quad
Alexander Yu. Olshanskii$^{\sharp\flat\natural}$
}
\address{
$^\flat$%
Faculty of Mechanics and Mathematics of Moscow State University
\\
Moscow 119991, Leninskie gory, MSU.
\\
$^\natural$Moscow Center for Fundamental and Applied Mathematics.
\\
$^\sharp$Department of Mathematics, Vanderbilt University, Nashville
37240, U.S.A.
\\
klyachko@mech.math.msu.su
\vphantom{$\sum^{\sum}$}
\qquad
werunik179@gmail.com
\qquad
alexander.olshanskiy@vanderbilt.edu
}
\grantsFirst{\RSF 22-11-00075}
\grantsThird{NSF, grant no.~DMS-1901976}
\abstract{
We show, in particular,
that, if a finite group $H$
is a retract of any
finite group
containing $H$ as a verbally closed subgroup,
then the centre of $H$ is a direct factor of $H$.
}

\s 1.
Introduction

A subgroup $H$ of a group $G$ is called \emph{verbally closed}
[MR14]
if
any equation of the form
$$
w(x,y,\dots)=h,
\qbox{where $w$ is an element of the free group
$F(x,y,\dots)$ and $h\in H$,
}
$$
having solutions in $G$ has a solution in $H$.
If each finite system of equations with coefficients from~$H$
$$
\{w_1(x,y,\dots)=1, \dots, w_m(x,y,\dots)=1\},
\qbox{where $w_i\in H*F(x,y,\dots)$ (and $*$ means the
free product),}
$$
having solutions in $G$ has a solution in $H$, then the subgroup $H$
is called \emph{algebraically closed} in $G$.

Algebraic closedness is a stronger property than verbal
closedness;
however these properties turn out to be
equivalent in many cases (see
[Rom12],
[RKh13],
[MR14],
[Mazh17],
[RKhK17],
[KM18],
[KMM18],
[Mazh18],
[Bog18],
[Bog19],
[Mazh19],
[RT19],
[RT20],
[Tim21]).
A group $H$ is called \emph{strongly verbally closed}
[Mazh18]
if it is algebraically closed
in any group containing $H$ as a verbally closed
subgroup. Thus, the verbal closedness is a subgroup property,
while the strong verbal closedness is a property of an abstract
group.
The class of strongly verbally closed groups is fairly wide.
For example, the following groups are strongly verbally closed:
\-
all abelian groups [Mazh18],
\-
all free groups [KM18],
\-
all virtually free groups
containing no
nonidentity finite normal subgroups
[KM18], [KMM18],
\-
all groups
decomposing nontrivially into a free product [Mazh19],
\-
fundamental groups of all connected surfaces
except the Klein bottle [Mazh18], [K21],

\enditem
See also [Bog18] and [Bog19]
for generalisations of
some of these results.
In particular,
these results imply that
\disp{\sl
the infinite dihedral group
is strongly verbally closed.
}%
This innocent looking
special case is one of the most difficult to prove
and requires a special argument.
Actually, this particular proposition is
the main result of
[KMM18] (used afterwards in [Mazh19],
[Bog18], and [Bog19]).
In Section~5, we complement
this fact with a description
of finite dihedral strongly verbally closed groups.

Proving the strong verbal closedness of a group is not easy,
but disproving this property is
not easy either. The literature contains only
the following examples of non-strongly-verbally-closed groups:
\-
the already mentioned fundamental group of the Klein bottle [K21]
\-
and two nonabelian groups of order eight [KM18], [RKhK17].

\enditem
We generalise the latter example substantially:
in Section 3, we prove that
\disp{\sl
the centre of any
finite strongly verbally closed group is its
direct factor
}%
(and, more generally,
strong verbal closedness of a finite group
implies stringent
constraints on its
abelian normal subgroups, see the centre theorem in Section~3).
In particular, this means that finite nilpotent
nonabelian groups cannot be strongly verbally closed.
This leads us to the following question.

\Question 1.
Does there exist a
finitely generated
nilpotent nonabelian strongly verbally closed
group?

We conjecture that the answer is negative;
but so far, we can prove this
only for the simplest infinite nilpotent
nonabelian group (see Section 4):
\disp{\sl
the Heisenberg group
$
\UT_3(\Z)=\pmatrix{
1&\Z&\Z
\cr
0&1&\Z
\cr
0&0&1
\cr
}
$
is not strongly verbally closed.
}

The algebraic closedness can be
characterised structurally if the group $H$ is
\emph{equationally Noetherian} (i.e. any system of
equations over $H$ with finitely many unknowns
has the same solutions in $H$ as some its
finite subsystem);
namely,
the algebraic closedness in this case is equivalent  to
the ``local retractness" [KMM18]:
\disp{
\sl\narrower
an equationally Noetherian
subgroup $H$ of a group $G$ is algebraically closed in $G$ if and only if
$H$ is a retract {\rm \(i.e. the image of an endomorphism $\rho$
of $G$
such that
$\rho\o\rho=\rho$\)} of each finitely generated over $H$
subgroup of $G$ \rm(i.e. a subgroup of the form $\gp{H\cup X}$,
where the set $X\subseteq G$ is finite).
}%
Another well-known fact
was proven in [RKh13] (Lemma 1.1):
\disp{\sl\hfuzz50pt 
if $V(G)$ is a verbal subgroup of a group $G$ and $H$ is a verbally
closed subgroup of~$G$, then
$H\cap V(G)=V(H)$
\(i.e. the verbal subgroup of~$H$
corresponding to the same variety\)
and
$H/V(H)\subseteq G/V(G)$
is verbally closed in $G/V(G)$.
}%
These two facts
imply
that
\disp{\sl\narrower
a finite group $H$ is strongly
verbally closed if and only if it is a retract of any finite group $G$
containing $H$ as a verbally closed subgroup
and satisfying all identities of $H$
}%
(because the variety generated by a finite group
consists of locally finite groups [Neu69], Theorem 15.71).
The following section contains
a (short)
proof that a finite group with nonabelian monolith
(e.g., any finite simple group) is
not only strongly verbally closed but also has a stronger
property.

\Th 1.
Any finite group with nonabelian monolith is strongly verbally
closed. Moreover, each such group $H$ is a retract of any
finite group containing $H$ and satisfying all identities of~$H$.
\newline
{\rm (So, the verbal-closedness requirement
can be omitted --- it holds
automatically in this situation.)}

\noindent
Recall that the \emph{monolith} of a group is the intersection of all
nonidentity normal subgroups of this group.
Many finite groups with abelian monoliths are also strongly verbally
closed and even retracts in groups from the corresponding varieties
(see Theorem 2 in the next section);
in particular,
semidirect products corresponding to
irreducible faithful representations
of finite $p'$-groups over $\Z_p$ enjoy this property.


We conclude this short introduction
with the following simple observation
(which, in particular, explains
the words ``finitely generated" in the open question above).

\proclaim{Embedding theorem}.
Any group
$H$ embeds into a strongly verbally closed group of
cardinality ~$|H|+\aleph_0$ that
satisfies all identities of $H$.

\Proof
If a group
$H$
is algebraically closed in a variety
$\cal M$
(i. e.
algebraically closed in any group from
$\cal M$ that contains $H$),
then $H$
is strongly verbally closed by
the lemma from [RKh13] mentioned above.
Thus,
the embedding theorem
follows from Scott's theorem [Sco51]:

\disp{\sl
any group $H$ from any variety $\cal M$ embeds into an
algebraically closed in $\cal M$ group from $\cal M$
of cardinality $|H|+\aleph_0$.
}%
In [Sco51]
(see also [LS80])
this theorem is stated for the variety of all groups,
but it is easy to see that
the
(simple)
proof
works without any modifications
for any variety
(and even for any class of groups
closed with respect to unions of ascending chains of groups).

\medskip

The authors thank Filipp Denissov for valuable remarks.

\s 2.
Monolithic groups

Let us recall some terminology (see [Neu69]):
\-
the \emph{variety generated by a class of groups $\cal K$}
is the class of all
groups satisfying all identities that hold in all groups from $\cal K$;
\-
a group $X$ is called \emph{critical} if it is not contained in the
variety generated by the class consisting of all proper
subgroups of $X$ and all quotient groups of $X$ by nontrivial normal
subgroups.

\enditem
\noindent{\bf Proof of the nonabelian-monolith theorem.}
Let
$M$ be the nonabelian monolith
of a
subgroup $H$ of a finite group $G$
satisfying all identities of~$H$.
We should show that $H$ is a retract
of $G$.
Any nontrivial normal subgroup
in a minimal hypothetical counterexample~$G$
nontrivially intersects $H$ (and, hence, $M$).
This implies that $G$ has a unique minimal normal
subgroup, and it contains $M$.
Both groups $H$ and $G$ are critical, because
their monoliths are nonabelian ([KoN66], see also [Neu69], 53.44).
Since
they generate the same variety, they are isomorphic
([Neu69], 53.33), i.e. $G=H$. This completes the proof.

For groups with abelian monoliths, the situation becomes more complicated.

\Lemma 1.
Let $M\ne\0$ be a finite module over
a group ring $\Z_{p^m}[G]$,
where the prime $p$ does not divide $|G|$. Then
\item{\rm 1)}
if the exponent of $M$ is $p^k$,
and all proper submodules has smaller exponents, then
\itemitem{\rm a)}
$M$ is homogeneous \(i.e., a direct sum of cyclics
of exponent $p^k$\),
\itemitem{\rm b)}
all its layers $p^iM/p^{i+1}M$, where $i\in\{0,\dots,k-1\}$,
are
isomorphic simple $G$-modules;
\itemitem{\rm c)}
elements of $G$ acting identically on the lower layer
act identically on the entire module $M$;
\item{\rm 2)}
if $M$ does not decompose into a direct sum nontrivially,
then all its proper submodules have smaller exponents.

\Proof
\item{\rm 1)}
{%
Put $L=\{x\in M\;|\;p^{k-1}x=0\}$.
By Maschke's theorem,
we have a decomposition
$M/pM=L/pM\oplus N/pM$
for a $G$-module $N$.
The submodule $N$ is of
exponent $p^k$, because
$L$ has smaller exponent.
Therefore, $L=pM$ by virtue of the minimality condition 1).
This proves~a).
\newline
The module $M/pM$ is simple
because of the same minimality condition:
the preimage in $M$
of any nonzero submodule of $M/pM$
has exponent $p^k$.
By virtue of the homogeneity,
the mapping
$$
M/pM\to M/p^{s+1}M,
\quad
a+pM\mapsto p^sa+p^{s+1}M
$$
is a well-defined isomorphism.
This proves b).
The property c) follows from b):
an element acting identically on the lower layer
acts identically on all layers;
and an automorphisms of an abelian $p$-group
acting identically on layers
form a $p$-group;
therefore, c) follows from the condition that
$|G|$
is not a multiple of $p$.
}

\item{\rm 2)}
Let us choose a minimal submodule $N\subseteq M$ of exponent $p^k$
equal to the exponent of $M$.
The homogeneity implies that $N$ has a direct complement in
the abelian
group $M$. Hence, it has a direct module complement
(this is a simple
generalisation of Maschke's theorem, see, e.g., [Pas83], Lemma 1.1).
Since $M$ is indecomposable, we obtain that $N=M$.

\Theorem 2.
Suppose that
a finite group
$H$ contains a normal subgroup $C$ such that
$C$
coincides with its centraliser, does not decomposes
into a direct
product of
nontrivial
subgroups
normal in $H$,
and $|C|$ is coprime to
$|H/C|$.
Then $H$ is a retract of any finite
group $G\supseteq H$ satisfying all identities of $H$.
In particular, $H$ is strongly verbally closed.

\Proof
Lemma 1 shows that $H$ is monolithic,
$C$ is the centraliser of the monolith
$M=\{c\in C\;|\;c^p=1\}$
and decomposes into a direct product
of cyclics of equal order $p^k$.
Note also that the subgroup $C$ is Sylow,
because its order is coprime to the index.

We can assume that
any nontrivial normal subgroup
of $G$ contains $M$, and $G$ is
monolithic with a monolith $L\supseteq M$
(because this is valid for any minimal hypothetical counterexample
$G$).

Note that $M$ is contained in the abelian
group $S=H^m$, where $m=|G:C|$.
Since $G\in\var N$, we have that
$M\subseteq G^m=S$,
where $S$
is a (normal) abelian Sylow subgroup of $G$, and,
hence, $L\subseteq S$. The abelian Sylow subgroup $S$ is a direct factor
of its centraliser by the Schur--Zassenhaus theorem:  $\CC_G(S)=S\times
U$. The subgroup $U$ is characteristic (and even verbal) in $\CC_G(S)$,
therefore, $U\nin G$. But $U\cap M=\1$,
i.e. $U=\1$, since $G$ is monolithic.
Thus, the Sylow $p$-subgroup~$S$ of $G$
coincides with its centraliser.
Therefore, Lemma 1 implies that $S$ is a
homogeneous group of
exponent $l$, and $l=k$, because $G\in\var N$;
moreover,
$S=\CC_G(L)$ (by Assertion 1)~c) of Lemma~1).

To prove that $G=H$,
it remains to show that
$|L|=|M|$ and $|G/\CC_G(L)|=|H/\CC_H(M)|$.
Since the left-hand sides of
these equalities do not exceed the right-hand sides
(because $H$ is a subgroup of $G$),
the equalities follow from
Lemma 53.25 of [Neu69], which
says, in particular, that
\disp{\sl
the monolith $L$ of any group $G$ from the variety generated by
a finite group $H$ is isomorphic to the monolith $N$ of some
factor~$F$
(i.e. a quotient group of a subgroup) of $H$,
and $G/\CC_G(L)\iso F/\CC_F(N)$.
}%
This means that, in the case under consideration, $L=M$
and $G/\CC_G(L)\iso H/\CC_H(M)$, which completes the proof.




Theorems 1 and 2
suggest the following definition.
We call a group
$H$
a
\emph{strong retract}
if it is a retract
of any
group $G\supseteq H$ from the variety $\var H$.
The property of being a
strong retract is a property of an abstract group $H$
(stronger than the strong verbal closedness).

\Proposition.
A finite subgroup $H$
of a group $G$
is a retract if and only if $H$ is a retract of any
finitely generated
subgroup of $G$
containing $H$.
In particular, the following groups are strong retracts:
\-
every finite group with nonabelian monolith
\-
and any finite group  $H$
containing a normal subgroup $C$ such that
$C$
coincides with its centraliser, does not decomposes
into a direct
product of
nontrivial
subgroups
normal in $H$,
and $|C|$ is coprime to
$|H/C|$.

\Proof
The assertion ``In particular" follows immediately from the main assertion
and Theorems 1 and 2.
The main assertion
is a corollary of Mal'cev's local theorem
(see [KaM82], Theorem 24.3.1):
\disp{\sl\narrower
if an object-universal formula $\Phi$ is true on subsystems
which form a local covering of an algebraic system
$A$, then $\Phi$ is true on $A$.
}%
It suffices to apply this theorem to the formula
``there exists a retraction $G\to H$" on the
algebraic system
$(G,\;\cdot,\;^{-1},\;h_1,\dots,h_n)$,
where $\{h_1,\dots,h_n\}=H$;
note that a retraction $\rho\:G\to H$
corresponds to a binary predicate
$P(x,y)$ (the graph of the retraction),
and the condition ``$P$ is a graph of a retraction"
can be written as a universal formula
(because $H$ is finite).

\Question 2.
What is an arbitrary finite strong retract?


\s 3.
The centres of finite strongly verbally closed groups

\proclaim Approximation lemma.
For any finite elementary abelian $p$-group $C$
\(where $p$ is prime\)
and any positive integer~$k$,
there exists $t\ge k$ such
that
the direct product $P=\bigtimes\limits_{i=1}^t C_i$
of copies $C_i$ of
$C$
contains a subgroup~$R$
invariant with respect to
the diagonal
action on $P$ of the
endomorphism algebra
$\End C\iso\MM_d(\Z_p)$
\(where $d=\log_p|C|$\) of $C$
with the following properties:
\item{\rm a)}
$R$ is contained in the union of kernels $K_j$ of the natural
retractions (projections) $P\to C_j$,
\item{\rm b)}
but $R\cdot\bigtimes\limits_{j\notin J}C_j=P$ for any subset
$J\subseteq\{1,\dots,t\}$ of cardinality $k$;
\item{\rm b$'$)}
moreover, each such $J$ is contained in a set $J'\supseteq J$
such that $P=R\times\bigtimes\limits_{j\notin J'}C_j$;
and there exist integers $n_{ij}$ such that the
projection
$\pi\:P\to\bigtimes\limits_{j\notin J'}C_j$
with kernel $R$
acts as:
$C_i\ni c_i\mapsto
\prod\limits_j c_j^{n_{ij}},
$
where $c_j\in C_j$ is the element corresponding to
$c_i$ under the isomorphism $C_i\iso C\iso C_j$.

\Proof
We regard the group $C$ (assumed to be nontrivial)
as the additive group of the
vector space $\Z_p^d$,
and the group~$P$ --- as the group of all mappings
from $X=V\setminus\0$ to $\Z_p^d$, where $V$ is a finite vector
space (of dimension $n$ chosen below)
over $\Z_p$ (so, $C_i$ is the set of
mappings vanishing at all but one vectors of
$X$, i.e. $t=p^n-1$).
As $R\subseteq P$, we take the set of
mappings defined by polynomials of
degree at most $r$
without free terms:
$$
\eqalign{
&R=R_r=
\cr
&=
\Bigl\{
(\^x_1,\dots,\^x_n)
\mapsto
\bigl(f_1(\^x_1,\dots,\^x_n),\dots,f_d(\^x_1,\dots,\^x_n)\bigr)
\;|\;f_i\in\Z_p[x_1,\dots,x_n],\;\deg f_i\le r,\;f_i(0,\dots,0)=0
\Bigr\}.
}
$$
Here $(\^x_1,\dots,\^x_n)$ are coordinates of a vector $v\in V$
with respect to some basis; but it is easy to see that the definition of
$R$ does not depend either on the choice of a basis in $V$
or on the choice of a basis in~$\Z_p^d$
that ensures invariance of $R$ with respect to
the ``diagonal"
action of $\GL_d(\Z_p)$.
By the same reason, $R$ is $\MM_d(\Z_p)$-invariant
(for any field $F$,
$\GL_d(F)$-invariance
implies $\MM_d(F)$-invariance,
because any matrix decomposes into a sum of nonsingular matrices).

The Chevalley theorem ([Che35], see also,
e.g., [Lan68], Chapter~5, Exercise~6) says that
\disp{\hfuzz255pt
\sl
polynomials
$f_1,\dots,f_d\in F[x_1,\dots,x_n]$ without free terms
over a finite field $F$ has a common nonzero root
if $n>\sum\deg f_j$.
}%
If we choose the space $V$ (i.e. the integer $t=p^{\dim V}-1=p^n-1$)
and the integer $r$ such that
$n=\dim V>rd$, then
property a) holds,
because any mapping $f\in R$ vanishes at a vector
from $X$, i.e. the projection of $f$ to a factor $C_i$
is
zero.

To prove b), we should find,
for any vectors
$v_1,\dots,v_k\in X$ and $w_1,\dots,w_k\in\Z_p^d$,
a mapping
$$
f=f_{v_1,w_1,\dots,v_k,w_k}\in R
$$
such that $f(v_i)=w_i$ for all $i$
(because then, for any mapping~$g\in P$,
we have $g=f_{v_1,g(v_1),\dots,v_k,g(v_k)}+h$,
where $h(v_i)=0$, as required).

Since the definition of $R$ does not depend on the choice of basis in $V$,
we can assume
that, for all vectors~$v_i$, all coordinates,
but the first
$k$ ones, are zero.
Any mapping from a $k$-dimensional vector space
over~$\Z_p$ to~$\Z_p$ is defined by a polynomial
of
degree at most $p-1$ in each variable, i.e. the total degree of such
polynomial is at most~$k(p-1)$.
Therefore, there exist polynomials~$f_1,\dots,f_d\in\Z_p[x_1,\dots,x_k]$
such that, for all $i\in\{1,\dots,k\}$ and $j\in\{1,\dots,d\}$,
$$
\eqalign{
(\hbox{the value of
$f_j$ at the first $k$ coordinates of $v_i$})
&=
(\hbox{$j$th coordinate of $w_i$}),\cr
f_j(0,\dots,0)
&=0,
\qqbox{and}
\deg f_j\le k(p-1).
}
$$
Therefore, if one choose
$r\ge k(p-1)$,
then the mapping $f\:(\^x_1,\dots,\^x_n)\mapsto
\bigl(f_1(\^x_1,\dots,\^x_k),\dots,f_d(\^x_1,\dots,\^x_k)\bigr)
$
lies in $R$ and maps $v_i$ to $w_i$,
which completes the proof of b).

Condition b$'$) follows from b),
because
the modules
$C_j$ are irreducible submodules of the semisimple
$\MM_d(\Z_p)$-module~$P$; therefore, it remains to
recall the following
well known (and easy-to-prove) general fact about modules:
\disp{\sl\hfuzz16pt
if $R$ is a submodule of a semisimple
module $P=R+\bigoplus\limits_{i\in I} C_i$,
where the modules $C_i$ are irreducible,
then $P=R\oplus\bigoplus\limits_{i\in I'} C_i$
for some subset $I'\subseteq I$.
}%
To explain
that the projection of $P$ onto $\bigtimes\limits_{j\notin J'}C_j$
is
defined by an integer matrix,
it
suffices to note that this projection
is an endomorphism of the $\MM_t(\Z_p)$-module and apply the
following well-known general fact on endomorphisms of semisimple
modules:
\disp{\sl
any endomorphism of
a module $\bigoplus C_i$ over a semisimple algebra $A$,
where
modules~$C_i$ are irreducible,
is defined by a matrix with entries from the centre of
$A$,
}%
i.e., for any endomorphism $\phi$, there exists a matrix
$M_\phi=(m_{ij})$ with elements from the centre of $A$ such
that, for each $c_i\in C_i$, we have
$\phi(c_i)=\displaystyle\sum\limits_{j:\;C_j\iso C_i} m_{ij}c_j$,
where $c_j\in C_j$ are elements corresponding
to $c_i$ under a (fixed) isomorphism~$C_j\iso C_i$.
This completes the proof of the lemma.

A similar proposition holds for any finite fields.
We
state it for a future reference (but
do not use it in this paper).
The proof can be easily obtained from the argument above
by obvious
alterations.

\proclaim Approximation lemma for finite fields.
For
any finite vector space $C$
over a finite field $F$
and any positive integer~$k$,
there exists $t\ge k$ such
that
the direct sum $P=\bigoplus\limits_{i=1}^t C_i$
of copies $C_i$ of
$C$
contains a subspace $R$
invariant with respect to
the
diagonal
action on $P$ of the
endomorphism algebra
$\End C\iso\MM_d(F)$
\(where $d=\dim C$\) of $C$
with the following properties:
\item{\rm a)}
$R$ is contained in the union of kernels $K_j$
of the natural
projections $P\to C_j$,
\item{\rm b)}
but $R+\bigoplus\limits_{j\notin J}C_j=P$ for any subset
$J\subseteq\{1,\dots,t\}$ of cardinality $k$;
\item{\rm b$'$)}
moreover, each such $J$ is contained in a set $J'\supseteq J$
such that $P=R\oplus\bigoplus\limits_{j\notin J'}C_j$;
and there exist $n_{ij}\in F$ such that
the
projection
$\pi\:P\to\bigoplus\limits_{j\notin J'}C_j$
with kernel $R$
acts as
$C_i\ni c_i\mapsto
\sum\limits_j n_{ij}c_j,
$
where $c_j\in C_j$ are vectors corresponding to~$c_i$ under the
isomorphism $C_i\iso C\iso C_j$.

\proclaim Centre theorem.
The centre of any finite strongly verbally closed group $H$
is a
direct factor of $H$.
Moreover, for any normal subgroup $N$
of a strongly verbally closed
group~$H$,
the centre~$\ZZ\Bigl(\CC\bigl(\ZZ(N)\bigr)\Bigr)$
of the centraliser~$\CC\bigl(\ZZ(N)\bigr)$
of the centre $\ZZ(N)$ of $N$
is a direct factor of this centraliser,
and some complement is
normal in~$H$\:
$$
\underbrace{\CC\bigl(\ZZ(N)\bigr)}%
_{\displaystyle{\cup\vrule height1.4ex \atop N}}
=
\underbrace{\ZZ\Bigl(\CC\bigl(\ZZ(N)\bigr)\Bigr)}%
_{\displaystyle{\cup\vrule height1.4ex \atop\ZZ(N)}}%
\times D
\qbox{for a subgroup $D\nin H$.}
$$

\Proof
The vertical inclusions are valid for any subgroup of any group
(we added them to the statement for clarity).

Denote $L=\CC\bigl(\ZZ(N)\bigr)$.
It suffices,
for each prime $p$, find a
homomorphism
$\psi_p\:L\to\ZZ(L)$,
commuting with the action of~$H$ on $L$
by
conjugations
and
injective on the $p$-component
$\ZZ_p(L)$ of the
centre $\ZZ(L)$ of $L$,
because then the homomorphism
$\psi\:x\mapsto\prod\limits_p\psi_p\bigl(\pi_p(x)\bigr)$
(where $\pi_p\:\ZZ(L)\to\ZZ_p(L)$ is the projection onto the $p$-component)
is injective on~$\ZZ(L)$ and, therefore, its kernel is the required
complement $D$.

Suppose that there is no such
homomorphisms $\psi_p$ for some $p$, i.e. any
homomorphism
$L\to\ZZ(L)$
commuting with the action of~$H$
is not injective on $\ZZ_p(L)$ and, hence, on
the maximal elementary subgroup~$C\subseteq\ZZ_p(L)$.
We should show that the group $H$ is not strongly verbally closed.

Let us choose $t$ by the approximation lemma
applied to~$C$
(for some integer $k$ to be specified later),
and
consider the fibered product
$$
Q=\left\{(h_1,\dots,h_t)\in
H^t\;\Bigm|\;h_1L=\dots=h_tL\right\}
\qbox{of $t$ copies of $H$.}
$$
As an overgroup $G\supset H$, we take the
quotient group~$G=Q/R$, where the subgroup $R\subset C^t$
is
chosen by the approximation lemma
($R$ is normal in $Q$, because $R$ is invariant with respect to
the diagonal
action of~$\Aut C$; and the conjugation action of $Q$ on
$P$ is diagonal, because $L^t$ commutes with $P=C^t$).

Then $H$ embeds into~$G$ diagonally: $h\mapsto(h,\dots,h)$.
This homomorphism is
an embedding not only into~$Q$ but also into~$G$
by the property a) of the approximation lemma, because
all projections of a nontrivial diagonal element of $Q$ are nontrivial.

This diagonal subgroup $H\subset G$ is
verbally closed in $G$.
Indeed, if an
equation $w(x_1,\dots,x_n)=h$ is solvable in $G$ and
$(\~x_1,\dots,\~x_n)\in H^t$ is a preimage of
a solution, then $w(\~x_1,\dots,\~x_n)=(hc_1,\dots,hc_t)$,
where $(c_1,\dots,c_t)\in R$ and, therefore, by the property a)
$c_i=1$ for some $i$, i.e.
the $i$th coordinates of the tuple~$(\~x_1,\dots,\~x_n)$ form
a solution to the
equation $w(x_1,\dots,x_n)=h$ in $H$, as required.

It remains to show, that the diagonal subgroup $H\subset G$
is not a retract. Let $\rho\:G\to H$ be a
hypothetical retraction and let $\^\rho\:Q\to H$ be
its composition with the natural epimorphism $Q\to Q/R=G$.

\noindent
Henceforth,
the
words ``subgroup" and ``centraliser"
refer to the fibered product $Q$
(which contains the diagonally embedded subgroup $H$);
the centraliser of a set~$X$ in $H$ is denoted by
$\CC_H(X)$, i.e. $\CC_H(X)=\CC(X)\cap H$.
If $U$ is a subgroup of~$L$,
then the symbol $U_i$, where $i=1,\dots,t$, denotes
the corresponding
subgroup
$\{(1,\dots,1,u,1,\dots,1)\;|\;u\in U\}$
of~$Q$.

Let us verify that
$$
\^\rho(L_i)\subseteq\CC_H\bigl(\CC_H(L)\bigr)=L
\qbox{for each $i$}.
\eqno{(*)}
$$
Indeed,
if $h\in C_H(L)$, then $h$ commutes with each
component of each element of $L$, i.e. $h$ commutes with~$L_i$.
Applying the retraction $\^\rho$ to this relation, we obtain
that $h=\^\rho(h)$ commutes with $\^\rho(L_i)$, which
proves the inclusion in~$(*)$.
To prove the equality, note that
$L$ is a centraliser (of a centre of $N$),
and the triple centraliser of any subgroup in any group coincides with
the single centraliser.

On the other hand, the
mutual commutator subgroup
$[L_i,\;L_j]$ is trivial for $i\ne j$.
Therefore, we obtain
$
[\^\rho(L_i),\;\^\rho(L_j)]=\1.
$
Hence, we have
$
\[\^\rho(L_i),\;\prod\limits_{j\ne i}\^\rho(L_j)\]=\1.
$
If $\^\rho(L_i)=\^\rho(L_l)$ for some different $i$ and~$l$,
then $\[\^\rho(L_i),\;\prod\limits_j\^\rho(L_j)\]=\1$ and, therefore,
$[\^\rho(L_i),\;L]=\1$ (because
$L=\^\rho(L)\subseteq\prod\limits_j\^\rho(L_j)$).
Thus, we obtain
$$
\^\rho(L_i)\subseteq\CC_H(L),
\qbox{if $\^\rho(L_i)=\^\rho(L_l)$ for some different $i$ and $l$}.
\eqno{(**)}
$$
Formulae $(*)$ and $(**)$ implies immediately that
\disp{\sl
if $\^\rho(L_i)=\^\rho(L_l)$ for some different $i$ and $l$,
then $\^\rho(L_i)\subseteq\ZZ(L)$.
}%
Let us take
$k$ in the approximation lemma to be the number of all subgroups of $H$,
and let $J$ be the set of \emph{exclusive}
numbers
$i$,
i.e. such
that
$\^\rho(L_i)\ne\^\rho(L_l)$
for any $l\ne i$.
Then, according to b$'$), we have a decomposition
$$
\bigtimes_{i=1}^t C_i
=
R\times\bigtimes\limits_{i\in I}C_i,
\qbox{where all $i\in I$ are non-exclusive},
$$
and the projection
$\pi\:\bigtimes\limits_{i=1}^t C_i\to\bigtimes\limits_{i\in I}C_i$
onto the second factor of this decomposition is defined by an integer matrix
$(n_{ij})$ (i.e.
$C_i\ni c_i\mapstoo^\pi \prod\limits_j c_j^{n_{ij}}$,
where $c_j\in C_j$ are elements corresponding to
$c_i$ under the isomorphism $C_i\iso C\iso C_j$).
This means that
the restriction $\^\pi\:C\to\bigtimes\limits_{i\in I}C_i$
of $\pi$ to $C$ is defined by the formula
$$
\^\pi\:c\mapstoo^\pi\prod\limits_{j\in I} c_j^{m_j},
\qbox{where $m_j={\sum\limits_i n_{ij}}$ and
$c_j\in L_j$ are
elements corresponding to $c\in C$.}
$$
Therefore, the composition
$$
\Psi\:C\to\ZZ(L),
\quad
\Psi\:
c\mapstoo^{\pi}\prod\limits_{j\in I} c_j^{m_j}
\mapstoo^{\^\rho}
\prod\limits_{j\in I}\^\rho\(c_j^{m_j}\)
$$
extends to a homomorphism $\Phi\:L\to\ZZ(L)$
defined by the similar formula:
$$
L\ni g\mapstoo^\Phi
\prod\limits_{j\in I}\^\rho\(g_j^{m_j}\),
\qbox{where $g_j\in L_j$ are
elements corresponding to $g\in L$.}
$$
(This is a homomorphism, because $\^\rho(L_j)$,
for $j\in I$,
is contained in the abelian
group $\ZZ(L)$.)
The homomorphism $\Phi$ commutes with the action of $H$
and, therefore, its kernel intersects $C$ nontrivially by
the assumption.
Thus, $\Psi$, being the restriction of
$\Phi$ to $C$, has a nontrivial kernel.
On the other hand, $\Psi$ is just the identity mapping,
as can be seen from its definition
($\^\rho\o\pi=\^\rho$, because $\^\rho(R)=\1$).
The obtained contradiction completes the proof.

\s 4.
Dihedral groups

We shall use the following simple facts:
\-
{\sl
any group $G$ from the variety generated by a dihedral group
of finite order
$2n$ or $4n$, where $n$ is odd,
decomposes into a semidirect product
$C\semitimes Q$, where $C$ is an elementary abelian 2-group
(Sylow 2-subgroup of~$G$), and $Q$ is an abelian group of exponent $n$
(Hall $2'$-subgroup of~$G$ or, equivalently, the verbal
subgroup generated by squares of all elements of~$G$);
thus, $Q$ is a $C$-module;
}
\-
{\sl
any finite module $V$ of
odd cardinality
over any elementary abelian 2-group $C$
decomposes into the direct sum
$V=\bigoplus\limits_{\chi\in\Chi}V_\chi$,
where $\Chi$ is the set of all homomorphisms
(\emph{characters}) $C\to\{\pm1\}$, and
$
V_\chi=\{v\in V\;|\;cv=\chi(c)v \hbox{ for all } c\in C\}
$.}

\proclaim Semidirect-product lemma.
In a semidirect product
$G=C\semitimes Q$, where $C$ is an elementary abelian 2-group
and $Q$ is a finite abelian group of odd exponent $n'$,
the equality
$$
\{g^{2n'/d}\;|\;g\in G\}=\{g\in G\;|\ g^d=1\}
\eqno{({**}*)}
$$
holds
for any divisor $d$ of $2n'$ such
that either $d=2$ or $d|n'$ and $\GCD(d,\;n'/d)=1$.

\Proof
Suppose that $g=cq$, where $c\in C$ and $q\in Q$.
According to the facts stated above we obtain the decomposition
$q=q_{c,+}q_{c,-}$, where
$cqc^{-1}=q_{c,+}q_{c,-}^{-1}$,
i.e. $q_{c,+}\in\prod\limits_{\chi\atop\chi(c)=1}Q_\chi$
and $q_{c,-}\in\prod\limits_{\chi\atop\chi(c)=-1}Q_\chi$.
Then
$$
g^k=(cq)^k=(cq_{c,+}q_{c,-})^k=
c^kq_{c,+}^kq_{c,-}^{{1\over2}\bigl(1+(-1)^k\bigr)}.
$$
This means that, for $d=2$,
the sets in both sides of $({**}*)$
consist of all possible products $cq_{c,-}$.
If $d$ is odd, then these sets are
$\{q^{n'/d}\;|\;q\in Q\}$ and $\{q\in Q\;|\ q^d=1\}$.
In an abelian group~$Q$ of exponent $n'$, these sets
coincide (if $d$ and $n'/d$ are coprime):
the both sets are the first direct factor in the
decomposition~$Q=\{q\in Q\;|\ q^d=1\}\times\{q\in Q\;|\;q^{n'/d}=1\}$.
This completes the proof.

\proclaim Dihedral-group theorem.
The dihedral group $D_n$ of order $2n$
is
strongly verbally closed
if and only if $n$ is either infinite
or not divisible by four.
\newline
Moreover,
suppose that $H=D_n=\gp b_2\semitimes \gp a_n$,
where $n$ is not a multiple of four, is a
dihedral subgroup
of a finite
group~$G$
belonging to the variety generated by
$H$.
Then the following conditions are equivalent:
\item{\rm1)}
$H$ is verbally closed in $G$;
\item{\rm2)}
$H$ is algebraically closed in $G$;
\item{\rm3)}
$H$ is a retract of $G$;
\item{\rm4)}
in terms of the decompositions from the beginning of the section
\($G=C\semitimes Q$ and $Q=\!\!\bigtimes\limits_\chi\! Q_\chi$\),
the order of the $\chi$-component $(a^2)_\chi$ of $a^2\in H$
equals the order of $a^2$
for some character
$\chi\:C\to\{\pm1\}$.

\Proof
The strong verbal closedness of infinite dihedral group is proven in
[KMM18].
Finite dihedral
groups~$D_{4k}$ are not strongly verbally closed
by the centre theorem.

It remains to prove the equivalence of Conditions 1) --- 4)
for finite $n$ not divisible by four.
To prove the implication $4)\imp3)$,
note that
the cyclic
(normal) subgroup $\gp{(a^2)_\chi}$
is a direct factor of $Q_\chi$
(because its order equals the period
of the abelian group $Q_\chi$)
and,
therefore, it is a direct factor of $Q$,
and the complement is normal in $G$
(because the action of~$G$ on~$Q_\chi$ is ``scalar" and,
hence,
all subgroups contained in $Q_\chi$ are normal in~$G$).
Let $\pi$ be the composition~$Q\to\gp{(a^2)_\chi}\too^\iso\gp{a^2}$
and let us
verify that the
mapping~$\phi\:G=C\semitimes Q\to H,\
c\cdot q\mapsto b^{{1\over2}(1-\chi(c))}\cdot\pi(q)$
is a homomorphism.
A
mapping from one semidirect product to another
$X\semitimes Y\to Z\semitimes T$
of the
form $xy\mapsto\alpha(x)\beta(y)$,
where $\alpha\:X\to Z$ and $\beta\:Y\to T$ are homomorphisms,
is a homomorphism if
$\beta(xyx^{-1})=\alpha(x)\beta(y)\alpha(x)^{-1}$
for all $x\in X$ and $y\in Y$.
In the case under consideration, this property holds:
$$
\pi(cqc^{-1})
=
\pi(cq_\chi c^{-1})
=
\pi(q_\chi^{\chi(c)})
=
\pi(q_\chi)^{\chi(c)}
=
b^{{1\over2}(1-\chi(c))}\pi(q_\chi)b^{-{1\over2}(1-\chi(c))}
=
\alpha(c)\pi(q_\chi)\alpha(c)^{-1}.
$$
\-
If $n$ is odd, then
$\phi$ is injective on $H$,
because any nontrivial normal subgroup of $H$
intersects $\gp a=\gp{a^2}$ nontrivially,
and the restriction of $\phi$ to $\gp{a^2}$ is injective
by Condition~4).
\-
If $n$ is even, then take any
homomorphism $\gamma\:G\to G/Q\to\gp{a^{n\over2}}$ such that
$\gamma(a^{n\over2})=a^{n\over2}$ and
consider the map
$\phi'\:G\to H=\gp{a^2,b}\times\gp{a^{n\over2}}$,
$g\mapsto\phi(g)\gamma(g)$, whose restriction to $H$
is obviously injective.

\enditem
Thus,
the composition of $\phi$ or $\phi'$
(depending on the parity of $n$)
with an automorphism of $H$
is the required retraction $G=C\semitimes Q\to H$.

The implications
$3)\imp2)\imp1)$ are general facts valid for any groups (see
Introduction).

\medskip

\noindent
It remains to prove the implication $1)\imp4)$,
i.e. to
construct an equation of the
form $w(x,y,\dots)=h$,
where~$w$ is an element of the free group
$F(x,y,\dots)$ and $h\in H$,
with the following properties:
\item{a)}
it is solvable in $G$;
\item{b)}
the solvability of this equation in $H$
implies Condition 4).

\enditem
In an explicit form, such an equation
looks fairly terrible.
For simplification, we do three things:
\-
First, we fix an element
$a'\in\gp a_n$ whose order is the product
of all odd prime divisors of $n$.
\-
Secondly, we shall construct a \emph{multi-sort equation}
(or an \emph{equation with typed variables}),
i.e. an equation in which a type $[d]$ is assigned
to each variable $x^{[d]}$,
\newline
where
$
d\in\{2\}\cup\{positive\ divisors\ d\ of\ n
\ coprime\ to\ n/d\}
$;
a solution
of a multi-sort equation
in a group
is a
substitution of elements of the group
(instead of the variables),
transforming the equation into a valid equality,
where a variable of type $[d]$ is allowed to be substituted only
by elements
of order dividing
$d$.
Let us verify that a multi-sort equation with
properties a) and b) can be
transformed to a usual equation with these properties;
indeed, if
we
replace each
typed variable~$x^{[d]}$ with
$x^{2n'/d}$,
where $n'=n$ if $n$ is odd,
and $n'=n/2$ if $n$ is even,
then we obtain a usual equation
whose
solvability (in~$G$ or $H$) is
equivalent to the solvability of the initial multi-sort equation
in the same group,
because $\{g^{2n'/d}\;|\;g\in G\}=\{g\in G\;|\ g^d=1\}$ by the
semidirect-product lemma.
\-
Thirdly, we write the equation in the module language,
i.e. the multiplication of elements from~$Q$ is denoted by $+$
(and the same symbol denotes the addition in the group ring of~$C$),
the symbol $\cdot$ denotes the
conjugation of elements of $Q$ by elements of $C$
from the left
(and also the multiplication in the group ring of $C$).

\enditem
For example, for $n=15$, the multi-sort equation
$
(x^{[2]}+2y^{[2]}x^{[2]})\cdot(3z^{[15]}+4t^{[5]})=a'
$
is translated into the usual language as
$
x^{15}(z^{6}t^{24})x^{-15}
(yx)^{15}(z^{6}t^{24})^2(yx)^{-15}
=a.
$
Certainly, the inverse translation is not always possible;
but we do not need it:
the only multi-sort equation we use
will be written in the module language and
it has the form similar to the example above:
a sum of integer polynomials in variables of type $[2]$
multiplied by variables of odd types equals an
element of $Q\cap H$ (namely, a power of $a'$).

We need the following simple identity
(almost copied from [KMM18]):
$$
\(\prod_{c\in C}\Bigl(1+\chi(c)c\Bigr)\)\cdot q=
2^{|C|}q_\chi
\qbox{for any
character $\chi$ and any
$q\in Q$}.
\eqno{(1)}
$$
To prove it, note that the
$\chi$-component of the element $q$
in the left-hand side is
multiplied by
two $|C|$ times;
while the other components vanish,
because, for each character~$\chi'\ne\chi$, there exists
$c\in C$ such that $\chi(c)=-\chi'(c)$.

Let us decompose the elementary abelian group
$C$ into a direct product of cyclics
$
C=\gp{c_1}_2\times\dots\times\gp{c_m}_2;
$
for each $c=\prod c_i^{\epsilon_i}\in C$ (where $\epsilon_i\in\{0,1\}$),
we put
$
f_c(x_1^{[2]},\dots,x_m^{[2]})
=
\prod\(x_i^{[2]}\)^{\epsilon_i}
$
and consider
the multi-sort equation
$$
\sum_{\chi\in X}
\(\prod_{c\in C}\Bigl(1+\chi(c)f_c(x_1^{[2]},\dots,x_m^{[2]})\Bigr)\)
\cdot
y_\chi^{\bigl[l\left|\gp{a'_\chi}\right|\bigr]}
=
2^{|C|}a',
\eqno{(2)}
$$
where $l$ is the minimal positive integer such
that $l\left|\gp{a'_\chi}\right|$ is coprime
with~$n/(l\left|\gp{a'_\chi}\right|)$.
This equation is solvable in $G$: just take
$x_i^{[2]}=c_i$ and
$y_\chi^{\bigl[l\left|\gp{a'_\chi}\right|\bigr]}=a'_\chi$,
and Equation (2) becomes the
sum over all characters~$\chi$ of Identities (1)
with $q=a'_\chi$. This proves Property~a).

It remains to study Equation (2) in the dihedral group $H$.
A substitution
$x_i^{[2]}\to h_i\in H$ determines a homomorphism $\phi\:C\to H/\gp a$,
$c\mapsto f_c(h_1,\dots,h_m)\gp a$
and
a character $\chi\:C\too^\phi H/\gp a\too^\iso\{\pm1\}$.
Any term
in the left-hand side of (2) corresponding to any other character
$\chi'\ne\chi$ vanishes, because there exists $c\in C$ such that
$\chi'(c)\ne\chi(c)$ and, therefore, the element
$
1+\chi'(c)f_c(h_1,\dots,h_m)=1+\chi'(c)\phi(c)\in\Z[H/\gp a]
$
annihilate the $(H/\gp a)$-module $\gp a$. Therefore,
this substitution transforms (2) into
$$
ky_\chi^{\bigl[l\left|\gp{a'_\chi}\right|\bigr]}=2^{|C|}a'
\qbox{(where $k=2^{|C|}$, as it is easy to see).}
$$
If
$\left|\gp{a'_\chi}\right|<\left|\gp{a'}\right|$,
then
\-
an odd prime divisor $p$ of $n$ does
not divide $\left|\gp{a'_\chi}\right|$
(by the definition of $a'$);
\-
hence,
$p$ does not divide
$l\left|\gp{a'_\chi}\right|$
(by the definition of $l$),
\-
therefore,
$y_\chi^{\bigl[l\left|\gp{a'_\chi}\right|\bigr]}$
cannot be substituted by elements of orders divisible by $p$
(by the definition of the solvability of a multi-sort equation);
\-
so,
the equation is not solvable in the dihedral group
(because the order of the right-hand side is divisible by~$p$).

\enditem
Thus, if Equation (2) is solvable in the dihedral group, then
$\left|\gp{a'_\chi}\right|=\left|\gp{a'}\right|$.
But $\gp{a'}$ is the
socle (i.e. product of all minimal normal subgroups)
of the group $\gp{a^2}$, i.e. the injectivity
of the restriction of the homomorphism
$a^2\mapsto(a^2)_\chi$ to $\gp{a'}$
implies the injectivity of this homomorphism on
the entire group $\gp{a^2}$.
Therefore, the solvability of Equation (2) in the dihedral group
implies the equality
$
\left|\gp{(a^2)_\chi}\right|
=
\left|\gp{a^2}\right|,
$
as required.

\medskip

The following
simplest nontrivial
example
illustrates the proof
of the implication 1) $\imp$ 4).

\Example.
Suppose that
$G=
D_3\times D_5=
\bigl(\gp{b_3}_2\semitimes\gp{a_3}_3\bigr)\times
\bigl(\gp{b_5}_2\semitimes\gp{a_5}_5\bigr)
$
and $H=D_{15}=\gp{b}_2\semitimes\gp{a}_{15}$
embeds into $G$ diagonally: $b=b_3b_5$ and $a=a_3a_5$.
The subgroup $H$ is
not a retract of $G$;
so, there must exist an equation solvable in $G$,
but not in $H$.
Our argument allows us to construct such an
equation explicitly.

In the case under consideration,
$C=\{1,b_3,b_5,b_3b_5\}$ (the Klein four-group),
$
Q=\gp a_{15}=\gp{a_3}_3\times\gp{a_5}_5,
$
and $a'=a$.
There are four characters $C\to\{\pm1\}$:
\-
two ``important": $\tau$ and $\pi$:\quad
$\tau(b_3)=-1=-\tau(b_5)=\pi(b_5)=-\pi(b_3)$,
\-
and two ``unimportant":
the trivial character $\epsilon$ and
the character $\delta=\tau\pi$:\quad
$
\delta(b_3)=-1=\delta(b_5).
$

\enditem
Now,
$Q_\tau=\gp{a_3}$,
$Q_\pi=\gp{a_5}$,
$a_\tau=a_3$,
and
$a_\pi=a_5$
(and $Q_\epsilon=Q_\delta=\1$).
There are
\-
two variables of type [2]:
$x_3^{[2]}$ and $x_5^{[2]}$,
whose translation into the usual language
are
$x_3^{15}$ and $x_5^{15}$,
\-
and two odd-type variables:
$y_3^{[3]}$ and $y_5^{[5]}$,
whose translations are
$y_3^{10}$ and $y_5^6$.

\enditem
The words $f_c$ (where $c\in C$) are
$
f_1=1,
\quad
f_{b_3}=x_3^{15},
\quad
f_{b_5}=x_5^{15},
\qqbox{and}
f_{b_3b_5}=x_3^{15}x_5^{15}.
$
The expressions~$\bigl(1+\chi(c)f_c\bigr)\cdot q$ from~(2)
(where $q$ is an expression in variables taking value in $Q$)
are
translated from the module language to the group one as:
$$
qf_cq^{\chi(c)}f_c^{-1}
=\cases{
(qf_c)^2,&if $\chi(c)=+1$ (because $f_c^2=1$ always);
\cr
[q,f_c], &if $\chi(c)=-1$.
}
$$
Thus, Equation (2) takes the form
$$
\[\[\((y_3^{10})^2x_5^{15}\)^2,x_3^{15}\],x_3^{15}x_5^{15}\]
\cdot
\[\[\((y_5^{6})^2x_3^{15}\)^2,x_5^{15}\],x_3^{15}x_5^{15}\]
=a^{16}=a,
$$
where the factors in the left-hand side correspond to the characters $\tau$
and $\pi$ (the factors corresponding to unimportant characters
are omitted, because they are equal to 1 identically).
Using the distributivity, we simplify this equation slightly:
$$
\[[(y_3^{20}x_5^{15})^2,x_3^{15}]\cdot[y_5^{12}(x_3^{15})^2,x_5^{15}],
\;x_3^{15}x_5^{15}\]
=a.
$$
In $G$, there is a solution:
$x_3=c_3$,
$x_5=c_5$,
$y_3=a_3^2$,
$y_5=a_5$
(i.e. $y_3^{[3]}=a_3$ and $y_5^{[5]}=a_5$).
In the dihedral group $H$, there are no solutions, because
\-
if
$x_3^{15}\in\gp{a}\not\ni x_5^{15}$,
then the first inner commutator
$[(\dots)^2,x_3^{15}]$ is 1,
and the second inner commutator
$[(x_3^{15}y_5^{6})^2,x_5^{15}]$ lies in $\gp{a^3}$
(for any $y_i$);
i.e. the left-hand side of the equation lies in $\gp{a^3}$
and cannot be equal to the right-hand side;
\-
if
$x_3^{15}\notin\gp{a}\ni x_5^{15}$,
then the left-hand side lies in $\gp{a^5}$
(by similar reasons)
and again cannot be equal to the right-hand side;
\-
if
$x_3^{15}\in\gp{a}\ni x_5^{15}$,
then the left-hand side is 1
and cannot be equal to the right-hand side.
\-
if
$x_3^{15}\notin\gp{a}\not\ni x_5^{15}$,
then $x_3^{15}x_5^{15}\in\gp a$
and the left-hand side is 1
because of the outer commutator.

\enditem


\s 5.
The Heisenberg group

\proclaim Affine-bilinear-function lemma.
Suppose that $U\supseteq U'$ and $V$ are finitely generated free
abelian groups,
$f\:U\times V\to\Z$ is a bilinear function
whose restriction to $U'\times V$ has
rank at least two.
Then,
for any $u\in U$,
the
set
$
f(u+U',V)\subseteq\Z
$
is
a subgroup
containing $f(U',V)$.

\Proof
Since bases in $U'$ and $V$ can be chosen
independently from each other, we can assume
that the matrix of $f$ is diagonal, and each integer on the diagonal
divides the next,
because
every integer matrix can be transformed to such a form
(the \emph{Smith normal form})
by integer invertible elementary
transformations of rows and columns,
see, e.g.,
[Vin99].
Thus, we can assume that
the function
$
(u',v)\mapsto f(u+u',v)
$
in coordinates takes the form
$
(u',v)\mapsto\sum n_ix_i'y_i+\sum m_iy_i
$
(where $x_i'$ and $y_i$ are coordinates of vectors~$u'$~and~$v$,
respectively, $m_i,n_i\in\Z$, and $n_1|n_2|\dots$,
while every $m_i$ linearly depends on the vector $u$).

By Dirichlet's theorem,
if $n_i\ne 0$,
the arithmetic progression~$\{n_ix_i'+m_i\;|\;x_i'\in\Z\}$
contains a number of the form $p_i\cdot\GCD(n_i,m_i)$,
where $p_1$ is an arbitrarily large prime.
Thus, the primes $p_i$ can be chosen pairwise distinct.
Therefore, for suitable $x_i'$,
we have
$
\GCD(n_1x_1'+m_1,\;n_2x_2'+m_2,\dots)
=\GCD(n_1,\;m_1,\;m_2,\dots)
$
(because $n_1|n_2|\dots$ and $n_2\ne0$, as the rank of the form is at
least two).

Hence,
choosing appropriate $y_i$,
we obtain
$
(n_1x_1'+m_1)y_1+(n_2x_2'+m_2)y_2+\dots=
\GCD(n_1,\;m_1,\;m_2,\dots).
$
This shows that
$f(u+U',V)=
\GCD(n_1,\;m_1,\;m_2,\dots)\cdot\Z
\supseteq
n_1\cdot\Z=f(U',V)$,
as required.

\proclaim Heisenberg-group theorem.
The Heisenberg group $\UT_3(\Z)$ is
not strongly verbally closed.

\Proof
A \emph{verbal mapping} in a group $H$
is a mapping $H^s\to H$ of the form
$$
(h_1,\dots,h_s)\mapsto w(h_1,\dots,h_s),
$$
where
$w(t_1,\dots,t_s)$ is an element of the free group
$F_s=F(t_1,\dots,t_s)$.

Let us verify that
\dispno{\sl
in the Heisenberg group $H=\UT_3(\Z)$,
the intersection
of the image of each verbal mapping
with the centre is a subgroup;
and the image of this mapping is a union
of some cosets of this subgroup.
}(3)%
Consider a verbal mapping
$\phi\:(h_1,\dots,h_s)\mapsto w(h_1,\dots,h_s)$
and put
$$
h_i
=T(x_i,y_i,z_i)
\:=
\pmatrix{
1&x_i&z_i
\cr
0&1&y_i
\cr
0&0&1
\cr
}.
$$
Then
$$
\phi(h_1,\dots,h_s)=
T\bigl(l(x_1,\dots,x_s),\;l(y_1,\dots,y_s),\;
f(x_1,\dots,x_s;y_1,\dots,y_s)+l(z_1,\dots,z_s)\bigr),
$$
where $l\:\Z^s\to\Z$ is a linear function and
$f\:\Z^s\times\Z^s\to\Z$ is a bilinear function.
By an automorphism of the free group $F$,
any word $w$ can be reduced to a normal form
$w(t_1,\dots,t_s)=t_1^mw'(t_1,\dots,t_s)$,
where $w'$ lies in the commutator subgroup of $F$.
In the normal form, the function $l$ is just
$l(x_1,\dots,x_s)=mx_1$, and
the restriction of
$f$ to $U'\times V$
is skew-symmetric,
i. e.
$f(u',u')=0$
for all $u'\in U'$;
indeed, if $u'=(a_1,\dots a_s)$
(and, hence, $a_1=0$ or $m=0$),
then
$
f(u',u')
$
is determined by the equality
$$
w\bigl(T(a_1,a_1,0),\dots,T(a_s,a_s,0)\bigr)
=
w'\bigl(T(a_1,a_1,0),\dots,T(a_s,a_s,0)\bigr)
=
T\bigl(0,0,f(u',u')\bigr)
$$
and, therefore, is zero,
since matrices of the form $T(a,a,0)$ commute
and the word~$w'$ is commutator.

\noindent
To prove (3), we should show that
\item{a)}
$\phi(H^s)\cap\ZZ(H)$ is a subgroup
\item{b)}
and,
for each element
$h=T(a,b,c)
=
\phi(\^h_1,\dots\^h_s)
\in\phi(H^s)
$,
the coset $h\cdot\bigl(\phi(H^s)\cap\ZZ(H)\bigr)$
is contained in $\phi(H^s)$.

\enditem
The both facts follow from the affine-bilinear
function lemma
applied to the free abelian groups
$$
U=\Z^s\supseteq U'=V=\{u\in\Z^s\;|\;l(u)=0\}.
$$
The condition  $\rk f|_{U'\times V}\ge2$ holds, because
the rank of a skew-symmetric function is always even
(if the rank is zero, we have nothing to prove).

Let us verify a).
We have
$
\phi(H^s)\cap\ZZ(H)
=
\Bigl\{
T\bigl(0,\;0,\;f(u,v)+l(w)\bigr)\;\Bigm|\;
u,v,w\in\Z^s,\;l(u)=l(v)=0
\Bigr\}.
$
By the lemma,
$\{f(U',V)\}$ is a subgroup of $\Z$.
Hence, $\{f(U',V)\}+l(\Z^s)$ is a subgroup too
(because $l(\Z^s)$ is obviously a subgroup).

Let us verify b).
We can assume that $b=0$,
because any element $h=T(a,b,c)$ can be reduced to such a form
by an
automorphism of $H$
(since $H$ is a
free nilpotent of class two,
and
any pair of matrices $\{T(a_1,b_1,c_1),\;T(a_2,b_2,c_2)\}$,
with $a_1b_2-a_2b_1=\pm1$
forms a free basis).
%
Now we obtain b) as follows.
If $\^h_i=T(\^x_i,\^y_i,\^z_i)$,
$u=(\^x_1,\dots,\^x_s)$
and
$v=(\^y_1,\dots,\^y_s)$,
then
$$
h\cdot\bigl(\phi(H^s)\cap\ZZ(H)\bigr)
=
T\bigl(l(u),\;0,\;f(u,v)+f(U',V)+l(\Z^s)\bigr)
\subseteq
T\bigl(l(u),\;0,\;f(u+U',V)+l(\Z^s)\bigr)
$$
(where the inclusion follows from the lemma).
For matrices~$\~h_i=T(\^x_i+u'_i,v'_i,z_i')$
we have
$$
\phi(\~h_1,\dots,\~h_s)=
T\bigl(l(u),\;0,\;f(u+u',v')+l(q')\bigr),
$$
where
$u'=(u'_1,\dots,u'_s)\in U'=V$,
$v'=(v_1',\dots,v_s')\in U'=V$, and
$q'=(z_1',\dots,z_s')\in U$ are
vectors, which we can choose.
Thus,
$
\phi(H^s)
\supseteq
T\bigl(l(u),\;0,\;f(u+U',V)+l(\Z^s)\bigr)
\supseteq
h\cdot\bigl(\phi(H^s)\cap\ZZ(H)\bigr).
$
This completes the proof of (3).

\medskip

\noindent
Now, let us take the central product
$
G=H\mathop\times\limits_{\scriptscriptstyle\ZZ(H)=\ZZ(\~H)}\~H=
(H\times\~H)/\{c\~c^{-1}\;|\;c\in\ZZ(H)\}
$
of
$H$ and its copy~$\~H$.
The verbal closedness of $H$ in $G$ follows
immediately from (3):
if $w\bigl((h_1,h_1'),\dots,(h_s,h_s')\bigr)=(h,1)$
in $G$, then
$w(h_1,\dots,h_s)=hc^{-1}$ and $w(h_1',\dots,h_s')=c\in\ZZ(H)$
for some $c\in\ZZ(H)$;
and then, according to (3), $h$ lies also in the image of the corresponding
verbal mapping.
Thus, the verbal closedness is proven.

The group $H$ is equationally Noetherian, as well as any linear group
[BMR99]; therefore, if $H$ is strongly verbally closed,
then it is a retract of $G$ (see Introduction).
But the retraction $G\to H$ is impossible, because $\~H$ commutes
with~$H$, hence,  $\~H$ must be mapped in the centre
by a hypothetical retraction,
but then $\ZZ(\~H)=\ZZ(H)$ is mapped to~$\1$.
This contradiction completes the proof.

\References


[BMR99]
G. Baumslag, A. Myasnikov, V. Remeslennikov,
Algebraic geometry over groups I. Algebraic sets and ideal theory,
J. Algebra, 219:1 (1999), 16-79.

[Bog18]
O. Bogopolski,
Equations in acylindrically hyperbolic groups and verbal closedness,
Groups, Geometry, and Dynamics (to appear).
\arXiv 1805.08071

[Bog19]
O. Bogopolski,
On finite systems of equations in acylindrically hyperbolic groups,
\newline
arXiv:1903.10906.

[Che35]
C. Chevalley,
D\'emonstration d'une hypoth\`ese de M. Artin,
Abh. Math. Semin. Univ. Hambg., 11 (1935), 73-75.

[KaM82]
M. I. Kargapolov, Ju. I. Merzljakov,
Fundamentals of the theory of groups.
Graduate Texts in Mathematics, 62, Springer, 1979.

[K21]
A. A. Klyachko,
The Klein bottle group is not strongly verbally closed,
though awfully close to being so,
Canadian Mathematical Bulletin, 64:2 (2021), 491-497.
\arXiv 2006.15523


[KM18]
A. A. Klyachko, A. M. Mazhuga,
Verbally closed virtually free subgroups,
Sb. Math., 209:6 (2018), 850-856.
\arXiv 1702.07761

[KMM18]
A. A. Klyachko, A. M. Mazhuga, V. Yu. Miroshnichenko,
Virtually free finite-normal-subgroup-free groups
are strongly verbally closed,
J. Algebra, 510 (2018), 319-330.
\arXiv 1712.03406

[KoN66]
L. G. Kov\'acs, M. F. Newman,
On critical groups,
J. Austral. Math. Soc., 6:2 (1966), 237-250.



[Lan68]
S. Lang,
Algebra,
Addison-Wesley Publishing Co., Inc., Reading, Mass. 1965.


[LS80]
R. Lyndon, P. Schupp,
Combinatorial group theory,
Springer, 2015.

[Mazh17]
A. M. Mazhuga,
On free decompositions of verbally closed subgroups
of free products of finite groups,
J.~Group Theory, 20:5 (2017), 971-986.
\arXiv 1605.01766

[Mazh18]
A. M. Mazhuga,
Strongly verbally closed groups,
J. Algebra, 493 (2018), 171-184.
\arXiv 1707.02464


[Mazh19]
A. M. Mazhuga,
Free products of groups are strongly verbally closed,
Sb. Math., 210:10 (2019), 1456-1492.
\arXiv 1803.10634



[MR14]
A. Myasnikov, V. Roman'kov,
Verbally closed subgroups of free groups,
J. Group Theory, 17:1 (2014), 29-40.
\arXiv 1201.0497


[Neu69]
H. Neumann,
Varieties of groups,
Springer-Verlag, Berlin-Heidelberg-New York, 1967.



[Pas83]
D. S. Passman,
It's essentially Maschke's theorem,
The Rocky Mountain Journal of Mathematics, 13:1 (1983), 37-54.


[Rom12]
V. A. Roman'kov,
Equations over groups,
Groups - Complexity - Cryptology,
4:2 (2012), 191-239.


[RKh13]
V. A. Roman'kov, N. G. Khisamiev.
Verbally and existentially closed subgroups of free nilpotent groups,
Algebra and Logic, 52:4 (2013), 336-351.


[RKhK17]
V. A. Roman'kov, N. G. Khisamiev, A. A. Konyrkhanova,
Algebraically and verbally closed subgroups and retracts
of finitely generated nilpotent groups,
Siberian Math. J., 58:3 (2017), 536-545.


[RT19]
V. A. Roman'kov, E. I. Timoshenko,
On verbally closed subgroups of free solvable groups,
Vestnik Omskogo Universiteta, 24:1 (2019), 9-16
(in Russian).


[RT20]
V. A. Roman'kov, E. I. Timoshenko,
Verbally Closed Subgroups of Free Solvable Groups,
Algebra and Logic, 59:3 (2020), 253-265.
\arXiv 1906.11689



[Sco51]
W. R. Scott,
Algebraically closed groups,
Proc. Amer. Math. Soc., 2:1 (1951), 118-121.

[Tim21]
E. I. Timoshenko,
Retracts and verbally closed subgroups
with respect to relatively free soluble groups,
Sib. Math. J., 62:3, 537-544 (2021).


[Vin99]
E. B. Vinberg,
A course in algebra,
Graduate Studies in Math. 56, Amer. Math. Soc.,
Providence, RI, 2003.

\end